\documentclass[12pt,oneside]{amsart}
\usepackage{amsfonts}
\usepackage{amscd}
\usepackage{amsthm}
\usepackage{amssymb}
\usepackage{amsmath}
\usepackage{epsfig}
\usepackage{graphicx}

\title{Simplicial volume of closed locally symmetric spaces of non-compact type}
\author{Jean-Fran\c{c}ois Lafont}
\address{Department of Mathematical Science,
Binghamton University,
Binghamton, NY 13902-6000}
\email{jlafont@math.binghamton.edu}
\author{Benjamin Schmidt}
\address{Department of Mathematics,
University of Michigan,
2074 East Hall, 530 Church St.,
Ann Arbor, MI 48109-1043}
\email{bischmid@umich.edu}

\theoremstyle{proposition}
\newtheorem{Lem}{Lemma}[section]

\newtheorem*{Def}{Definition}

\theoremstyle{plain}

\newtheorem{Thm}[Lem]{Theorem}

\theoremstyle{remark}

\newtheorem*{Prf}{Proof}

\begin{document}

\begin{abstract}
We show that compact, locally symmetric spaces of non-compact type 
have positive simplicial volume.  This gives a positive answer to a question
that was first raised by Gromov \cite{Gr} in 1982.  We provide a summary
of results that are known to follow from positivity of the simplicial volume.
\end{abstract}

\maketitle

\section{Introduction}
In his paper \cite{Gr}, Gromov introduced the notion of the
simplicial volume of a connected, closed, and orientable
manifold $M$. This homotopy invariant is denoted by
$||\,M\,||\in [0,\infty)$ and measures how efficiently the
fundamental class of $M$ may be represented using real cycles.  


In the same paper, the question was raised as to whether the
simplical volume of a closed locally symmetric space of non-compact
type is positive (pg. 11 in \cite{Gr}).  Since then, this question
has been mentioned in a variety of different sources (\cite{Gr2},
\cite{L}, \cite{Sa}, \cite{CF1}), and has become a well-known
``folk conjecture.''  The purpose of this paper is 
to answer this conjecture in the affirmative.  Namely, we obtain:

\vskip 10pt

\noindent{\bf Main Theorem:} If $M^n$ is a closed locally symmetric
space of non-compact type, then $||M^n||>0$.

\vskip 10pt

The approach we use is due to Thurston \cite{Th} and
bounds the simplicial volume $||M||$ from below by the proportion
between the volume of $M$ and the maximal volume of suitably
defined straightened top dimensional simplices in $M$.
Technically, however, we are indebted to Besson,
Courtois, and Gallot for their pioneering work around the use of
the barycenter method in proving the rank one entropy rigidity
conjecture for locally symmetric spaces \cite{BCG} and to 
Connell and Farb for their subsequent development of the
technique in higher rank spaces (see \cite{CF1} for an extensive
survey).  

The main contribution of this paper lies in the idea of using
the barycenter method in order to define the straightened
simplices that are central to Thurston's argument. This allows us
to obtain control over the straightening process from estimates
(similar to those) in \cite{CF2}.  

We would like to thank Dick
Canary, Tom Farrell, and Ralf Spatzier for their interest in our
work and many helpful discussions.  We would also like to thank 
Wolfgang L\"uck for some helpful information concerning the application 
in section 4.11.

\section{Background}
In this section we collect together the relevant definitions and
results that are used in the proof of the main theorem.

\subsection{Simplicial volume}

We begin with the definition of the simplicial volume and the
important proportionality principle.

\begin{Def}
Let $M$ be a topological space, $C^0(\Delta^k, M)$ be the set of singular
$k$-simplices, and let $c=\sum_{i=1}^{j} a_i \cdot f_i$ with
each $a_i\in \mathbb{R}$ and $f_i \in C^0(\Delta^k, M)$ be a singular
real chain.  The $L^1$ norm of $c$ is defined by $||c||_{L^1} =
\sum_i |a_i|.$ The $L^1$ norm of a real singular homology class
$[\alpha]\in H_k^{sing}(M,\mathbb{R})$ is defined by $||\,
[\alpha]\,||_{L^1}=inf\{||c||_{L^1}\,|\, \partial(c)=0 \,,
\,[c]=[\alpha]\}$.
\end{Def}

\begin{Def}
Let $M^n$ be an oriented closed connected $n$-manifold with
fundamental class $[M^n]$. The simplicial volume of $M^n$ 
is defined as $||M^n||=||i([M^n])||_{L^1}$,
where $i:H_n(M, \mathbb{Z}) \rightarrow H_n(M,\mathbb{R})$ is the change of
coefficients homomorphism, and $[M^n]$ is the fundamental class
arising from the orientation of $M^n$.
\end{Def}

The proportionality principle (\cite{Gr},\cite{Th},\cite{St}) for simplicial
volume is expressed in the following:

\begin{Thm}
Let $M$ and $M'$ be two closed Riemannian manifolds with isometric
universal covers.  Then
$$\frac{||M||}{\textmd{Vol}(M)}=\frac{||M'||}{\textmd{Vol}(M')}.$$
\end{Thm}

In addition, the simplicial volume is particularly well behaved with respect
to products and connected sums.  Namely, the following relationships 
hold:

\begin{Thm}
For a pair of closed manifolds $M_1$, $M_2$, we have that:
$$C\cdot ||M_1||\cdot ||M_2||\geq ||M_1\times M_2|| \geq ||M_1||\cdot ||M_2||$$
where $C>1$ is a constant that depends only on the dimension
of $M_1\times M_2$.
\end{Thm}

\begin{Thm}
For $n\geq 3$, the connected sums of a pair of $n$-dimensional manifolds $M_1$
and $M_2$ satisfy:
$$||M_1\sharp M_2||=||M_1||+||M_2||$$
\end{Thm}

The proof of these last two results can be found in Sections 1.1 and 3.5 of \cite{Gr}
respectively (for Theorem 2.2, see also Section F.2 in \cite{BP}).  

\subsection{Thurston's approach}

Our proof will follow Thurston's method in \cite{Th} of showing positivity
of the simplicial volume of $M^n$.  Our version of Thurston's method can be summarized 
in the following:

\begin{Thm}
Let us denote by $\tilde M^n$ the universal cover of $M^n$, $\Gamma$ the
fundamental group of $M^n$, and
$C^0(\Delta^k,\tilde M^n)$ be the set of singular $k$-simplices in 
$\tilde M^n$, where $\Delta^k$ is assumed to be equipped with a fixed 
Riemmanian metric.  Assume that we are given a collection of maps 
$st_k:C^0(\Delta^k,\tilde 
M^n)\rightarrow C^0_{st}(\Delta^k,\tilde M^n)$ where $C^0_{st}
(\Delta^k,\tilde M^n)\subset C^0(\Delta^k,\tilde M^n)$.  We will
say this collection of maps is a {\it straightening} provided it satisfies 
the following four formal properties:
\begin{enumerate}
\item the maps $st_k$ are $\Gamma$-equivariant,
\item the maps $st_*$ induce a chain map $st_*:C^{sing}_*(\tilde M^n,\mathbb R)
\rightarrow C^{sing}_*(\tilde M^n, \mathbb R)$ which is $\Gamma$-equivariantly 
chain homotopic to the identity, 
\item the image of $st_n$ lies in $C^1(\Delta^n,\tilde M^n)$, i.e. 
straightened top-dimensional simplices are $C^1$,
\item there exists a constant $C>0$, depending solely on $\tilde M^n$ and the
chosen Riemannian metric on $\Delta^n$, such that 
for any $f\in C^0(\Delta^n,\tilde M^n)$, and corresponding straightened
simplex $st_n(f):\Delta^n\rightarrow \tilde M^n$, there is a uniform 
upper bound on the Jacobian of $st_n(f)$:
$$|Jac(st_n(f))(\sigma)|\leq C$$
where $\sigma \in \Delta^n$ is arbitrary, and the Jacobian is computed
relative to the fixed Riemannian metric on $\Delta^n$.
\end{enumerate}
If such a straightening exists, then $||M^n||>0$.
\end{Thm}

In the Theorem, one could replace properties (3) and (4) by the more general 
condition that the volume of the images of straightned top-dimensional simplices
are uniformly bounded above.  This more general approach was used in \cite{Th} 
and \cite{IY} to give a proof that closed negatively curved manifolds 
have positive simplicial volume. 
A different
straightening procedure was developed by Savage \cite{Sa} to show positivity
of the simplicial volume for cocompact quotients of $SL_n(\mathbb R)/SO_n(\mathbb R)$.
Our proof of the main theorem will involve a new straightening procedure
for locally symmetric spaces of non-compact type which do not have any local
$\mathbb H^2$ or $SL_3(\mathbb R)/SO_3(\mathbb R)$ factors.  Our formulation isolates
properties (3) and (4) because the barycenter method we will use is particularly
well-adapted to establishing these properties.

We finish with a brief explanation of how Theorem 2.4 is proved.  We first 
note that property (1) implies that the straightening procedure descends to
a straightening procedure on the compact quotient $M^n$.  Property (2) ensures
that the homology of $M^n$ obtained via the complex of straightened chains 
coincides with the ordinary singular homology of $M^n$.  Furthermore, since
the straightening procedure is a projection operator on the level of chains, 
it is contracting in the $L^1$-norm.  In particular, if $\sum a_if_i$ is a real 
$n$-chain representing the fundamental class of $M^n$, then so is $\sum a_ist(f_i)$,
and we have the inequality $||\sum a_if_i||_{L^1}\geq ||\sum a_ist(f_i)||_{L^1}$.

As a consequence, in order to show that the simplicial volume of $M^n$ is positive, 
it is sufficient
to give a lower bound for the $L^1$-norm of a straightened chain representing 
the fundamental class.  But now observe that, by property (3), the straightened
chain is $C^1$, and hence we can compute the volume of $M^n$ by:
$$Vol(M^n) = \int _{\sum a_ist(f_i)} dV_{M^n} =\sum a_i\int_{st(f_i)}dV_{M^n},$$
where $dV_{M^n}$ is the volume form on $M^n$. On the other hand, we have the bound:
$$\sum a_i\int_{st(f_i)}dV_{M^n}\leq \sum
|a_i|\int_{\Delta^n}|Jac(st(f_i))|dV_{\Delta^n},$$
where $dV_{\Delta ^n}$ is the volume form for the fixed Riemannian metric on $\Delta^n$.
Now by property (4) the Jacobian of straightened simplices is bounded uniformly 
from above, and hence we have a uniform upper bound: 
$$\int_{\Delta^n}|Jac(st(f_i))|dV_{\Delta^n}\leq K$$
where $K>0$ depends solely on $M^n$.  This yields the inequality:
$$Vol(M^n)\leq K\cdot \sum |a_i|$$
which upon dividing, and passing to the infimum over all straightened 
chains, provides the positive
lower bound $||M^n||\geq Vol(M^n)/K>0$.


\subsection{Locally symmetric spaces}

The principle technical tool in our argument is the use of the
barycenter method (see \cite{BCG} for a history as well as a survey of
applications). We will use this method in order to define our
straightenings. The technique ensures that straightened simplices
are $C^1$ and also gives a pointwise upper bound for their Jacobian
(formal properties (3) and (4) in Theorem 2.4). 


For details on locally symmetric spaces of non-compact type, we
refer the reader to the summary given in \cite{CF2}. Let
$(M^n,g_0)$ denote a closed locally symmetric space of non-compact
type. We let $(X,g)$ be the symmetric universal covering space of
$M$ and fix a basepoint $p\in X$. We will use the following notation:

\begin{itemize}

\item $G=\textmd{Isom}(X)^{0},$ $K=\textmd{Stab}_{G}(p)$, so that
$X\cong G/K$, and $P$ denotes a minimal parabolic subgroup of $G$.

\item The visual boundary and Furstenberg boundaries are denoted
by $\partial X$ and by $\partial_{F}X\cong G/P$, respectively.

\item The $h(g_0)$-conformal density given by the family of
Patterson-Sullivan measures is denoted by $ \nu : X \rightarrow
\mathcal{M}(\partial X)$, where $\mathcal{M}(\partial X)$ denotes 
the space of atomless probability measures on the visual boundary 
of $X$.

\item The map $\nu$ is $\Gamma$-equivariant, in the sense that $\nu(\gamma x)$
coincides with $\gamma_*\nu(x)$, the pushforward of the measure $\nu(x)$ under
the $\gamma$ action on $\partial X$.

\item $B_{p}:X \times
\partial X \rightarrow \mathbb{R}$ denotes
the Busemann function based at the point $p\in X$.

\end{itemize}

In higher
rank, the family of Patterson-Sullivan measures were constructed
in  \cite{K}, \cite{Al} and are $\pi_1(M)$-equivariant, atomless, and fully 
supported on $\partial_F X$. Additionally, the Busemann functions satisfy 
$B_{p}(\cdot,\cdot)=B_{\gamma p}(\gamma \cdot, \gamma \cdot)$ for each
 $\gamma \in \pi_{1}(M)$ and $p\in X$.  We will denote by:
$$dB_{(x,\theta)}(\cdot):T_x X \rightarrow \mathbb R$$
$$DdB_{(x,\theta)}(\cdot,\cdot):T_x X \otimes T_x X \rightarrow \mathbb R$$
the 1-form and 2-form obtained by differentiating the Busemann function (based
at the point $p\in X$)
corresponding to the direction $\theta\in \partial X$ at the point
$x\in X$.  Note that, while the Busemann {\it functions} depend on the 
chosen basepoint $p\in X$, the 1-form and 2-form defined above do {\it not}.
This justifies the omission of the basepoint $p\in X$ in our notation for
these forms.

For a measure $\mu \in\mathcal{M}(\partial X)$, let
$$g_{\mu}(\cdot)=\int_{\partial X} B_{p}(\cdot,\theta)\,
d(\mu)(\theta).$$  
When $g_\mu:X\rightarrow \mathbb{R}$ has a unique minimum, the
barycenter of $\mu$, denoted by $bar(\mu)\in X$, is defined to be the unique
point where $g_{\mu}$ is minimized.  In \cite{CF2}, the proof of
their proposition 3.1 shows that measures fully supported on the
Furstenberg boundary have well defined barycenters and that they are 
independent of the chosen basepoint $p \in X$ used to define $g_\mu$.

Finally, we recall the following result of Connell and Farb 
(see Section 4 in \cite{CF2}):

\begin{Thm}
Let $M$ be a closed locally symmetric space of non-compact type
with no local direct factors locally isometric to $\mathbb{H}^2$
or $SL_3(\mathbb{R})/SO_3(\mathbb{R})$, and let $X$ be its
universal cover. Let $\mu \in
\mathcal{M}(\partial X)$ be a probability measure fully supported
on $\partial_F X$ and let $x\in X$. Consider the endomorphisms $K_x(\mu)$, 
$H_x(\mu)$, defined on $T_xX$ by:
$$\langle K_x(\mu)u,u\rangle =\int_{\partial_{F}X}
DdB_{(x,\theta)}(u,u)\,d(\mu)(\theta)$$ and
$$\langle H_x(\mu)u,u\rangle =\int_{\partial_F X}
dB^2_{(x,\theta)}(u)\, d(\mu)(\theta).$$
Then $\det(K_x(\mu))> 0$ and there is a positive constant $C:=C(X)>0$
depending only on $X$ such that:
$$J_x(\mu):=\frac{\det(H_x(\mu))^{1/2}}{\det(K_x(\mu))}\leq C.$$
Furthermore, the constant $C$ is explicitly computable.
\end{Thm}

In the proof of our main theorem, this result of Connell and Farb
will be applied in the very specific case where $X$ is an irreducible
higher rank locally symmetric space, distinct from 
$SL_3(\mathbb R)/SO_3(\mathbb R)$, and the measure $\mu$ is a weighted
sum of Patterson-Sullivan measures.

\section{Proof of the Main Theorem}

We start out by providing a reduction, showing that it is sufficient to 
prove the Main Theorem for {\it irreducible} closed locally symmetric spaces 
of non-compact type.
In order to see this, let $M$ be an arbitrary locally symmetric
space of non-compact type, with universal cover $X$.  We observe that by the 
proportionality principle
in Theorem 2.1, in order to show that $||M||>0$, it is sufficient to show 
that $||M^\prime||>0$ for some locally symmetric space of non-compact
type whose universal cover is $X$.
 
Now let $G$ denote the identity component of $\textmd{Isom}(X)$
and $G=G_1\times\cdots\times G_k$ be the product decomposition of
$G$ into simple groups corresponding to the product decomposition
of $X$ into irreducible symmetric spaces. By a result of Borel \cite{Bo}, 
there are cocompact lattices $\Gamma_i \subset G_i$ for each $i\in
\{1,\ldots,k\}$.  Take $M^\prime$ to be the product locally symmetric 
space $M_1\times\cdots\times M_k$ obtained from the product lattice
$\Gamma_1 \times \cdots \times \Gamma_k$. 
From Theorem 2.2, the inequality $||M_1\times \cdots \times
M_k|| \geq \prod_{i=1}^{k}||M_i||$ holds.  Hence, if one has the 
Main Theorem for {\it irreducible} locally symmetric spaces of 
non-compact type, one obtains the Main Theorem for {\it all} locally 
symmetric spaces of non-compact type.

Next we observe that for the irreducible locally symmetric spaces 
modelled on $\mathbb{H}^n$, $\mathbb C\mathbb H^n$, $\mathbb H\mathbb H^n$,
and $Cay\mathbb H^2$ (the rank one cases), positivity of the simplicial
volume follows from \cite{Th} and \cite{IY}.  Furthermore, for 
$SL_3(\mathbb{R})/SO_3(\mathbb{R})$, 
positivity of the simplicial volume follows from \cite{Sa}.  
The Main Theorem will then follow from:

\vskip 5pt

\noindent{\bf Claim:} If the manifold $M^n$ is a compact quotient of an
irreducible
higher rank symmetric space of non-compact type $X$, and 
$X \neq SL_3(\mathbb{R})/SO_3(\mathbb{R})$, then the simplicial
volume satisfies $||M^n||>0$.

\vskip 5pt

To obtain the Claim, we use the Thurston approach.  From Theorem 2.4,
it is sufficient to define a straightening for $M^n$.  Before 
proceeding to do this, we fix some notation.
Recall that a singular $k$-simplex in $M^n$ is a continuous map
$f:\Delta^k \rightarrow M^n$, where $\Delta^k$ is the standard
Euclidean $k$-simplex realized as the convex hull of the standard
unit basis vectors in $\mathbb{R}^{k+1}$.  For our purpose it is
more convenient to work with the spherical $k$-simplex
$\Delta^k_s=\{ (a_1, \ldots,a_{k+1}) |\ a_i\geq 0,
\sum_{i=1}^{k+1} a_i^2=1\}\subset \mathbb{R}^{k+1}$, equipped with
the Riemannian metric induced from $\mathbb{R}^{k+1}$.  We will 
denote by $e_i$ ($1\leq i\leq k+1$) the standard basis vectors for
$\mathbb{R}^{k+1}$.  Finally, we will denote by $\Gamma$ the fundamental
group of $M^n$.

We now define the straightening procedure we will use:

\begin{Def} Given a singular $k$-simplex $f \in C^0(\Delta^k_s,X)$, 
with vertices $x_i:={f}(e_i)$, define $st_k ({f}) \in
C^{0}(\Delta^k_s,X)$ by $st_k ({f})(\sum_{i}a_i
e_i)= bar(\sum_{i} a_i^2 \nu(x_i))$. 
\end{Def}

The fact that the simplex $st_k (f)$ is well defined follows from the
comments in Section 2.3. Moreover, observe that $st_k (f)$ depends only
on the vertices of the original simplex $f$, i.e. on the $(k+1)$-tuple
of points $V:=(x_1,\ldots x_{k+1})$.  Let $V$ denote the collection of
vertices of $f$. As $st_k (f)$ depends only on $V$, we set
$st_V(\sigma):=st_k (f)(\sigma)$.  We now proceed to 
verify that this straightening procedure satisfies the four formal
properties needed.  For the convenience of the reader, we restate each 
property prior to proving it.

\vskip 10pt

\noindent{\bf Property (1):}
The maps $st_k$ are $\Gamma$-equivariant.

\begin{Prf}
Fix a point $\sigma=\sum_i a_i e_i \in 
\Delta^{k}_{s}$.  Then for any $\gamma \in \Gamma,$ $st_{\gamma V}(\sigma)$ 
is defined as the unique minimizer of the function $g_{\sigma}(\cdot)= 
\int_{\partial_F X} B_p (\cdot,\theta) d(\sum_i a_{i}^2 \nu(\gamma x_i))(\theta)$.  
Since 
$$ \int_{\partial_F X} B_p (\cdot,\theta) d(\sum_i a_{i}^2 \nu(\gamma x_i))(\theta)=$$ 
$$ \int_{\partial_F X} B_p (\cdot,\theta) d(\sum_i a_{i}^2 \gamma_{*}\nu(x_i))(\theta)=$$ 
$$ \int_{\partial_F X} B_p (\cdot,\gamma^{-1} \theta) d(\sum_i a_{i}^2 \nu(x_i))(\theta)=$$ 
$$ \int_{\partial_F X} B_{\gamma^{-1} \gamma p} (\gamma^{-1} \gamma \cdot,\gamma^{-1} 
\theta) d(\sum_i a_{i}^2 \nu(x_i))(\theta)=$$ 
$$ \int_{\partial_F X} B_{\gamma p} (\gamma \cdot,\theta) d(\sum_i a_{i}^2 
\nu(x_i))(\theta),$$ 
and since $B_{\gamma p}(\cdot,\cdot)$ and $B_p(\cdot,\cdot)$ differ by a function 
$k(\theta)$ of $\theta$, it follows that the unique minimizer of $g_{\sigma}(\cdot)$ 
is also the unique minimizer of the function:
$$h_{\sigma}(\cdot)= \int_{\partial_F X} B_{p} 
(\gamma \cdot,\theta) d(\sum_i a_{i}^2 \nu(x_i))(\theta).$$ 
Indeed, we have that the difference of the two functions is:
$$g_\sigma(\cdot)-h_\sigma(\cdot)=\int_{\partial_F X}k(\theta) d(\sum_i a_{i}^2 
\nu(x_i))(\theta)$$
which is a constant function on $X$.  But if $x\in X$ is the 
unique minimizer of $h_{\sigma}(\cdot)$, then $\gamma^{-1}x=st_V(\sigma)$.  It 
follows that $st_{\gamma V}(\sigma)=\gamma st_{V}(\sigma)$, completing the proof
of Property (1).
\end{Prf}

\vskip 10pt

\noindent{\bf Property (2):}
The maps $st_*$ induce a chain map $st_*:C^{sing}_*(X,\mathbb R)
\rightarrow C^{sing}_*(X, \mathbb R)$ which is $\Gamma$-equivariantly 
chain homotopic to the identity.

\begin{Prf}
The fact that $st_k$ commutes with the boundary operator follows
from the fact that $st_k(f)$ depends solely on the vertices of the
singular simplex $f$, along with the fact that $st_k(f)$ restricted
to a face of $\Delta^k_s$ coincides with the straightening of that 
face.

To see that $st$ is chain homotopic to the identity, first note
that the uniqueness of geodesics in $X$ gives rise to a well
defined $\Gamma$-equivariant straight line homotopy between any
simplex ${f}$ and its straightening $st({f})$. Hence
there are canonically defined homotopies between simplices and
their straightenings in $X$. Moreover, these homotopies when
restricted to lower dimensional faces agree with the homotopies
canonically defined on those faces. Appropriately ($\Gamma$-equivariantly)
subdividing these homotopies defines the required chain homotopy, 
concluding the proof of Property (2).
\end{Prf}

\vskip 10pt





\noindent{\bf Property (3):}
The image of $st_n$ lies in $C^1(\Delta^n_s,X)$, i.e. 
straightened top-dimensional simplices are $C^1$.

\begin{Prf}
Notice that for any simplex $f\in C^0(\Delta^n_s,X)$ and any $\sigma=\sum_i a_ie_i \in \Delta_s^n$, 
we have an implicit characterisation of the point
$st_n(f)(\sigma)=st_V(\sigma)$ via the 1-form equation:
\begin{equation} 0 \equiv
d(g_{\sigma})_{st_V(\sigma)}(\cdot) = \int_{\partial_{F}X}
dB_{(st_V(\sigma),\theta)}(\cdot) \, d(\sum_i a_i^2
\nu(x_i))(\theta). \end{equation}
Indeed, the fact that $st_V(\sigma)=bar(\sum_i
a_i^2 \nu(x_i))$ is defined as the unique minimum of the function:
$$g_{\sigma}(\cdot)=\int_{\partial_{F}X} B(p,\cdot,
\theta)\, d(\sum_i a_i^2 \nu(x_i))(\theta),$$
yields equation (1) upon differentiating.  

Since the map $st_V$ is given implicitly, one can apply the
implicit function theorem: in order for $st_V$ to be $C^1$, 
one needs to check the non-degeneracy condition.  But this 
merely requires that for the endomorphism $K$ defined by:
$$\langle K(u), u\rangle:=\int_{\partial_{F}X}
DdB_{(st_V(\sigma),\theta)}(u,u)\,d(\sum_i a_i^2\nu(x_i))(\theta)$$
defined on the tangent space $T_{st_V(\sigma)}M^n$, the determinant
be non-zero.  Note however
that in the notation of Theorem 2.5, the 
determinant of this matrix is precisely 
$\det(K_{st_V(\sigma)}(\sum_i a_i^2\nu(x_i)))$, and hence must be non-zero
as the measure $\sum_i a_i^2\nu(x_i)$ has full support on the Furstenberg
boundary.  This completes the proof of property (3).
\end{Prf}

\vskip 10pt

\noindent{\bf Property (4):}
There exists a constant $C>0$, depending solely on $X$, such that 
for any $f\in C^0(\Delta^n_s,X)$, and corresponding straightened
simplex $st_n(f):\Delta^n_s\rightarrow X$, there is a uniform 
upper bound on the Jacobian of $st_n(f)$:
$$|Jac(st_n(f))(\sigma)|\leq C$$
where $\sigma=\sum_i a_i e_i \in \Delta^n_s$ is arbitrary, and the Jacobian is computed
relative to the Riemannian metric on the spherical simplex $\Delta^n_s$
induced from $\mathbb R^{n+1}$.

\begin{Prf}
Differentiating the implicit 1-form equation with respect to directions in
$T_{\sigma}(\Delta_s^n)$, one obtains the two form equation
\begin{equation}
\begin{split}
0 \equiv D_{\sigma}d(g_{\sigma})_{st_V(\sigma)}(\cdot,\cdot) =
\sum_i 2a_i\langle\cdot,e_i\rangle_{\sigma}\int_{\partial_{F}X}
dB_{(st_V(\sigma),\theta)}(\cdot) \, d(\nu(x_i))(\theta) \\
+ \int_{\partial_{F}X}
DdB_{(st_V(\sigma),\theta)}(D(st_V)_{\sigma}(\cdot),\cdot)
\,d(\sum_i a_i^2 \nu(x_i))(\theta).
\end{split}
\end{equation}
defined on $T_\sigma(\Delta^n_s)\otimes T_{st_V(\sigma)}(X)$.  Now define symmetric 
endomorphisms $H_{\sigma}$ and $K_{\sigma}$ of
$T_{st_V(\sigma)}(X)$ by

$$\langle H_{\sigma}(u),u\rangle_{st_V(\sigma)}=\int_{\partial_{F}X}
dB^2_{(st_V(\sigma),\theta)}(u)\,d(\sum_i a_i^2\nu(x_i))(\theta)
\,\text{, and}$$
$$\langle K_{\sigma}(u),u\rangle_{st_V(\sigma)}=\int_{\partial_{F}X}
DdB_{(st_V(\sigma),\theta)}(u,u)\,d(\sum_i a_i^2\nu(x_i))(\theta).
$$

The fact that $K_{\sigma}$ is positive definite follows from
Theorem 2.5. Let $\{v_j\}_{j=1}^{n}$ be an orthonormal eigenbasis
of $T_{st_V(\sigma)}(X)$ for $H_\sigma$.  At points $\sigma \in
\Delta_s^n$ where the Jacobian of $st_V$ is nonzero, let
$\{\tilde{u}_j\}$ be the basis of $T_\sigma(\Delta_s^n)$ obtained
by pulling back the $\{v_j\}$ basis by $K_\sigma \circ D(st_V)_{\sigma}$, 
and $\{u_j\}$ be the orthonormal basis of
$T_\sigma(\Delta_s^n)$ obtained from the $\{\tilde{u}_j\}$ basis
by applying the Gram-Schmidt algorithm. We now have the sequence of
equations (which we will justify in the next paragraph):

\begin{equation}
\det(K_\sigma) \cdot |Jac(st_V)(\sigma)| =|\det(K_\sigma \circ D(st_V)_{\sigma})|
\end{equation}
\begin{equation}
=\prod_{j=1}^n |\langle K_\sigma \circ D(st_V)_{\sigma}(u_j),v_j\rangle_{st_V(\sigma)}|
\end{equation}
\begin{equation}
=\prod_{j=1}^{n}\Big|\sum_{i=1}^{n+1} \langle u_j,e_i\rangle_\sigma \cdot 2a_i
\int_{\partial_F X} dB_{(st_V(\sigma),\theta)}(v_j)\,
d(\nu(x_i))(\theta)\Big|
\end{equation}
\begin{equation}
\leq \prod_{j=1}^{n} \Big[\sum_{i=1}^{n+1}
\langle u_j,e_i\rangle_\sigma^2\Big]^{1/2}\,\Big[\sum_{i=1}^{n+1} 4a_i^2 
\Big(\int_{\partial_F X} dB_{(st_V(\sigma),\theta)}(v_j)\, d(\nu(x_i))(\theta)
\Big)^2\Big]^{1/2}
\end{equation}
\begin{equation}
\leq 2^{n} \prod_{j=1}^{n} \Big[\sum_{i=1}^{n+1}a_i^2 \int_{\partial_F X}
dB^2_{(st_V(\sigma),\theta)}(v_j)\, d(\nu(x_i))(\theta)\Big]^{1/2}
\end{equation}
\begin{equation}
=2^n\prod_{j=1}^{n}\langle H_\sigma (v_j), v_j\rangle_{st_V(\sigma)}^{1/2}
=2^{n}\det(H_\sigma)^{1/2},
\end{equation}

We now justify each step in the previous list of equations. 
Equation (3) follows from the definition of the Jacobian, along
with the fact that $\det(AB)=\det(A)\cdot \det(B)$.  Equation (4)
follows from the fact that, with respect to the
$\{u_j\}$ and $\{v_j\}$ bases, $K_\sigma \circ D(st_V)_{\sigma}$
is upper triangular, and hence the determinant is the product of the
diagonal entries.  Equation (5) follows from equations (4) and (2). 
Inequalities (6) and (7) follow from the Cauchy-Schwartz inequality applied in
$\mathbb{R}^{n+1}$ and the spaces $L^2(\partial_F X,\nu(x_i))$,
respectively, along with the fact that the $u_j$ are unit vectors
in $T_\sigma(\Delta^n_s)\subset T_\sigma(\mathbb R^{n+1})$.  The
two equalities in (8) follow from the definition of $H_\sigma$, and the fact
that the $\{v_j\}_{j=1}^{n}$ is an orthonormal eigenbasis for $H_\sigma$.

Upon dividing, we now obtain the inequality:
$$|Jac(st_V)(\sigma)|\leq 2^n\frac{\det(H_\sigma)^{1/2}}{\det(K_\sigma)}$$
But now note that, in the notation of Theorem 2.5, the expression 
$\det(H_\sigma)^{1/2}/{\det(K_\sigma)}$ is 
exactly $J_{st_V(\sigma)}(\sum a_i^2\nu(x_i))$.  Since the measure
$\sum a_i^2\nu(x_i)$ has full support in the Furstenberg boundary, Theorem 2.5 now yields
a uniform constant $C^\prime$, depending solely on $X$, with the property
that:
$$|Jac(st_V)(\sigma)|\leq 2^nJ_{st_V(\sigma)}(\sum a_i^2\nu(x_i))\leq 2^nC^\prime =: C$$
This completes the proof
of Property (4).
\end{Prf}



Having verified Properties (1)-(4) in the definition of straightening, 
we now conclude that the Claim holds, completing the proof of the Main Theorem.

\section{Applications}

In this section, we summarize the known consequences of 
positivity of simplicial volume.  
Most of the applications we mention can be found in Gromov's 
original paper \cite{Gr}.  We also refer the reader to Pansu's
article \cite{Pa} and to Chapter 14 in L\"uck's book \cite{L}.

For the applications we give, we point out that:
\begin{itemize}
\item those discussed in Sections 4.5-4.8, and 4.11 were
previously unknown for higher rank locally symmetric spaces of 
non-compact type.
\item the result in Section 4.1 was unknown for higher rank locally 
symmetric spaces of non-compact type that contain local $\mathbb H^2$ or
$SL_3(\mathbb R)/SO_3(\mathbb R)$ factors.
\item the estimates in Sections 4.9 and 4.10 can be explicitly 
computed, as our procedure gives a computable bound for the simplicial
volume.
\item let $\mathcal M$ be the smallest class of topological manifolds that 
(1) contains all closed locally symmetric spaces of non-compact type, (2)
is closed under connected sums with arbitrary closed manifolds of dimension
$\geq 3$, (3) is closed
under products, and (4) is closed under fiber extensions by surfaces of
genus $\geq 2$ (i.e. if $M\in \mathcal M$, and $M^\prime$ fibers over $M$
with fiber a surface $S_g$ of genus $\geq 2$, then $M^\prime \in \mathcal M$).
Then combining our Main Theorem, Theorems 2.2 and 2.3 from the introduction,
along with a result of Hoster and Kotschick \cite{HK}, one obtains that for
every manifold $M\in \mathcal M$, $||M||>0$.  For the manifolds in $\mathcal M$,
we obtain all the applications given in Sections 4.1 through 4.8.  We also
point out that by a result of Kapovich and Leeb \cite{KL}, there exist
surface bundles over surfaces (both of genus $\geq 2$, and hence in the
class $\mathcal M$) that do {\it not} support metrics of non-positive curvature.
More generally, manifolds in the class $\mathcal M$ that arise from a connected 
sum will fail to be
aspherical, and hence cannot support metrics of non-positive curvature.
\end{itemize}
We now list out the applications.

\subsection{Degree theorem}  

We provide a new application of positivity of the simplicial volume:

\begin{Lem} Suppose that $(N^n,g_{N})$ and
$(M^n,g_M)$ are connected, closed, and orientable Riemannian
manifolds and that $||M^n||>0$. If $||N^n||=0$, then there are no
continuous maps $f:N^n \rightarrow M^n$ of positive degree.
Otherwise, there is a constant
$C:=C(g_{\tilde{M}},g_{\tilde{N}})>0$ depending on the
Riemannian universal coverings of $M$ and $N$ such that
$$\textmd{deg}(f) \leq C
\frac{\textmd{Vol}(N)}{\textmd{Vol}(M)}.$$
\end{Lem}

\begin{Prf}
Assume $||M||>0$, and that $f:N\rightarrow M$ is a continuous map.  It
is easy to see that $||N|| \geq \deg(f) ||M||$ (pg. 8 in \cite{Gr}) and since
$||M||>0$, this gives the inequality $\deg(f)\leq ||N||/||M||$.  This immediately
implies that if $||N||=0$, then $\deg(f)=0$.

On the other hand, if $||N||>0$, then by Theorem 2.1, there are constants
$\delta_{g_{\tilde{M}}},\delta_{g_{\tilde{N}}}>0$ such that $||M||=\delta_{g_{\tilde{M}}}
\textmd{Vol}(M)$ and $||N||=\delta_{g_{\tilde{N}}}
\textmd{Vol}(N)$. Letting $C:=\delta_{g_{\tilde{N}}}/\delta_{g_{\tilde{M}}}$, we
immediately obtain:
$$\textmd{deg}(f)\leq \frac{||N||}{||M||}=
\frac{\delta_{g_{\tilde{N}}}\textmd{Vol}(N)}{\delta_{g_{\tilde{M}}}\textmd{Vol}(M)}=
C\frac{\textmd{Vol}(N)}{\textmd{Vol}(M)}$$
concluding the proof of the Lemma.
\end{Prf}

This yields a degree theorem (upper bound on the degree in terms of
the volume ratio between the codomain and domain) whenever the codomain manifold
has positive simplicial volume.  The question of 
whether the degree theorem could be obtained from positivity
of the simplicial volume was brought up in the survey paper
\cite{CF1}.

\subsection{Co-Hopf property}

The {\it co-Hopf property} for a group $G$ states that every monomorphism
$G\hookrightarrow G$ is in fact an isomorphism.  Note that
the co-Hopf property for $\pi_1(M^n)$ follows immediately from 
the Degree theorem, provided that $M^n$ is aspherical.  The co-Hopf property
for lattices was first shown by Prasad in \cite{Pr}, and also follows from 
Margulis' superrigidity theorem in the higher rank case.

\subsection{Positivity of MinVol}

The {\it minimal volume} of a smooth manifold $M$, denoted by 
$\textmd{MinVol}(M)\in [0,\infty)$, is defined as the infimum of
$\textmd{Vol}(M,g)$ as $g$ varies through complete Riemannian metrics with
$|\text{K}(g)|\leq 1$.  It was shown by Gromov (pgs. 35-37 in \cite{Gr})
that positive simplicial volume implies positive MinVol.  

\subsection{Positivity of Minimal Entropy}

The {\it minimal entropy} of a smooth manifold $M$ is defined to
be the infimum of the topological entropies of the geodesic flow
over all complete 
Riemannian metrics of unit volume on $M$.  There is 
the following inequality between simplicial volume and the minimal 
entropy $h$ (see pg. 37 in \cite{Gr}):
$$C\cdot ||M||\leq h(M)^n$$
where $C$ is a uniform constant, depending only on the dimension $n$
of $M$.  Hence positivity of simplicial norm implies positivity of
the minimal entropy. 

\subsection{Non-collapsing}

We say that $M$ {\it collapses} provided that there exists
a sequence of Riemannian metrics $g_i$ on $M$, satisfying $|\text{K}(g_i)|\leq 1$,
and having the property that at every point $p\in M$, the 
injectivity radius with respect to the metric $g_i$ is $<1/i$.
Gromov showed that manifolds with positive simplicial volume do {\it not}
collapse (pgs. 67-68 in \cite{Gr}).  

\subsection{Non-existence of F-structures}  

Loosely speaking, a positive rank {\it F-structure} on a manifold $M^n$ consists of a 
finite open cover $\mathcal U:= \{U_i\}$, along with effective torus actions $T^{k_i}$ 
(of dimension $k_i\geq 1$)
on each $U_i$, with the property that the torus actions commute on all 
the various intersections (we refer the reader to the survey article by
Fukaya \cite{F} for the precise definition as well as applications).  
A fundamental result of Cheeger and Gromov \cite{CG}
is that existence of positive rank F-structures is equivalent to collapsing.  Hence
if the simplicial volume of a manifold is positive, it does not support any
F-structure.  We point out that a smooth, fixed point free $S^1$-action is a 
special case of an F-structure. In particular, positivity of the 
simplicial volume implies that the manifold cannot support any such 
$S^1$-actions (see also the paper by Yano \cite{Y}). 

\subsection{Non-vanishing of top-dimensional bounded cohomology}

Bounded cohomology $\hat H^*(M^n)$ is defined in Gromov \cite{Gr}, where it 
is shown (pgs. 16-17) that $M^n$ has positive simplicial volume if and only if 
the map induced by inclusion of chains 
$i^n:\hat H^n(M^n)\rightarrow H^n_{sing}(M^n, \mathbb R)$ is non-zero.
This immediately implies that the $n$-dimensional bounded cohomology
of our manifolds $M^n$ is in fact non-zero. 

In fact, Gromov has shown (pgs. 46-47 in \cite{Gr}) the following
more general statement: if $f:X\rightarrow Y$ is
a continuous map between path-connected spaces, with the property that
the induced map $f_*:\pi_1(X)\rightarrow \pi_1(Y)$ is surjective
with amenable kernel, then the induced map $\hat f^*:\hat H^*(Y)
\rightarrow \hat H^*(X)$ on bounded cohomology is an isometric 
isomorphism.  In view of the previous paragraph, one obtains 
that the $n$-dimensional bounded cohomology $\hat H^n(X)$ is 
non-zero for any space $X$ which has a map $f:X\rightarrow M^n$
satisfying the hypotheses above.

\subsection{Covering by amenable subsets}

Let us call a subset $Y\subset M^n$ {\it amenable} provided that 
for every path component $Y^\prime\subset Y$, with inclusion map
$i:Y^\prime\rightarrow M^n$, the image of the morphism $i_*:\pi_1(Y^\prime)
\rightarrow \pi_1(M^n)$ is amenable.  It is shown by Gromov (pgs. 47-48 in 
\cite{Gr}) that if a space $X$ has a cover $\{U_\alpha\}$ by amenable sets,
with the property that every point $p\in X$ lies in $\leq k$ of the sets in
the cover, then the bounded cohomology $\hat H^i(X)=0$ for all $i\geq k$.

Hence the non-vanishing of the $n$-dimensional bounded cohomology of $M^n$
implies that for every covering $\{U_\alpha\}$ of $M^n$ by amenable sets,
there exists a point $p\in M^n$ contained in at least $(n+1)$ of the sets 
in our cover.  

\subsection{Bounds on Euler characteristics of flat bundles}

Let $E\rightarrow M^n$ be an $n$-dimensional affine flat bundle over $M^n$, and 
$\chi$ the Euler number of $E$ (obtained by evaluating the Euler class of $E$ on
the fundamental class of $M^n$).  Then the inequality $|\chi|\leq 2^{-n}||M^n||$
holds.  

This is due to an argument of Smillie, presented by Gromov (pgs. 21-23 in \cite{Gr}),
improving on the earlier inequality $|\chi|\leq ||M^n||$ due to Milnor and Sullivan
(\cite{Mi}, \cite{Su}).  

\subsection{Bounds on the sum of Betti numbers}

If $M$ is a connected sum of locally symmetric spaces of non-compact type,
then there exists a constant $C$, depending solely on the dimension of $M$,
with the property that:
$$\sum b_i(M) \leq C\cdot ||M||,$$
where the $b_i(M)$ are the Betti numbers of $M$ with any given coefficients.
This result is announced on pg. 12 of \cite{Gr}.

\subsection{$L^2$-invariants}

In L\"uck's book, the following question is raised (Conjecture 14.1, pg. 485, in
\cite{L}): let $M$ be an aspherical closed orientable manifold
of dimension $\geq 1$, and suppose that $||M||=0$.  Does it follow that 
$\tilde M$ is of determinant class, and satisfies $b_p^{(2)}(\tilde M)=0$ 
for all $p\geq 0$, and $\rho^{(2)}(\tilde M)=0$?

We observe that, in view of the computations on pg. 230 of \cite{L}, there
exist locally symmetric spaces of non-compact type that do {\it not} satisfy
the conclusion of this conjecture.  By our main theorem, they also fail to
satisfy the hypotheses of this conjecture, and hence do not provide 
counterexamples to this conjecture.

\section{Concluding remarks}

We conclude by pointing out some open questions related to our Main Theorem:

\vskip 5pt

\noindent{\bf Conjecture:} Let $M^n$ be a closed Riemannian manifold, whose
sectional curvatures are $\leq 0$, and whose Ricci curvatures are $<0$.  Then
$||M^n||>0$.

\vskip 5pt

This conjecture was attributed to Gromov in \cite{Sa}.  It seems plausible that a similar approach could be used to verify this 
conjecture.  The main difficulty lies in obtaining formal property (4) for the
analogous straightening procedure when the space $M^n$ is locally irreducible 
and is {\it not} a locally symmetric space.  We can also ask the:



\vskip 5pt





\noindent{\bf Question:} For a given closed locally symmetric space of non-compact
type $M^n$, what is the precise value of the ratio $||M^n||/\textmd{Vol}(M^n)$? 

\vskip 5pt

One of our applications is the non-vanishing of the top-dimensional 
bounded cohomology.  We have the natural:

\vskip 5pt

\noindent{\bf Question:} What is the dimension of $\hat H^n(M^n)$ for a locally 
symmetric space of non-compact type?  In particular, is it finite dimensional?

\end{document}